\documentclass[12pt]{amsart}
\usepackage{amsmath}
\usepackage{amsxtra}
\usepackage{amscd}
\usepackage{amsthm}
\usepackage{amsfonts}
\usepackage{amssymb}
\usepackage{eucal}
\usepackage{color}

\usepackage[all]{xy}

\newcommand{\nc}{\newcommand}
\newcommand{\rc}{\renewcommand}

\numberwithin{equation}{section}

\newtheorem{thm}{Theorem} [section]
\newtheorem{prop}[thm]{Proposition}

\newtheorem{conj}[thm]{Conjecture}

\newtheorem{important note}[thm]{Important Note}

\nc{\on}{\operatorname}

\nc{\bA}{{\mathbb A}}
\nc{\bB}{{\mathbb B}}
\nc{\bC}{{\mathbb C}}
\nc{\bD}{{\mathbb D}}
\nc{\bE}{{\mathbb E}}
\nc{\bF}{{\mathbb F}}
\nc{\bG}{{\mathbb G}}
\nc{\bH}{{\mathbb H}}
\nc{\bI}{{\mathbb I}}
\nc{\bJ}{{\mathbb J}}
\nc{\bK}{{\mathbb K}}
\nc{\bL}{{\mathbb L}}
\nc{\bM}{{\mathbb M}}
\nc{\bN}{{\mathbb N}}
\nc{\bO}{{\mathbb O}}
\nc{\bP}{{\mathbb P}}
\nc{\bQ}{{\mathbb Q}}
\nc{\bR}{{\mathbb R}}
\nc{\bS}{{\mathbb S}}
\nc{\bT}{{\mathbb T}}
\nc{\bU}{{\mathbb U}}
\nc{\bV}{{\mathbb V}}
\nc{\bW}{{\mathbb W}}
\nc{\bZ}{{\mathbb Z}}
\nc{\bX}{{\mathbb X}}
\nc{\bY}{{\mathbb Y}}

\newcommand{\bbC}{\mathbb{C}}

\newcommand{\bbZ}{\mathbb{Z}}
\nc{\Z}{\bbZ}
\nc{\C}{\bbC}

\nc{\cA}{{\mathcal A}}
\nc{\cB}{{\mathcal B}}
\nc{\cC}{{\mathcal C}}
\nc{\cD}{{\mathcal D}}
\nc{\cE}{{\mathcal E}}
\nc{\hE}{\hat{\cE}}
\nc{\hA}{\widehat{A}}
\nc{\hK}{\widehat{K}}
\nc{\hM}{\hat{\cM}}
\nc{\hN}{\hat{\cN}}

\nc{\cF}{{\mathcal F}}
\nc{\cG}{{\mathcal G}}
\nc{\cH}{{\mathcal H}}
\nc{\cI}{{\mathcal I}}
\nc{\cJ}{{\mathcal J}}
\nc{\cK}{{\mathcal K}}
\nc{\cL}{{\mathcal L}}
\nc{\cM}{{\mathcal M}}
\nc{\cN}{{\mathcal N}}
\nc{\cO}{{\mathcal O}}
\nc{\cP}{{\mathcal P}}
\nc{\cQ}{{\mathcal Q}}
\nc{\cR}{{\mathcal R}}
\nc{\cS}{{\mathcal S}}
\nc{\cT}{{\mathcal T}}
\nc{\cU}{{\mathcal U}}
\nc{\cV}{{\mathcal V}}
\nc{\cW}{{\mathcal W}}
\nc{\cZ}{{\mathcal Z}}
\nc{\cX}{{\mathcal X}}
\nc{\cY}{{\mathcal Y}}

\nc{\fA}{{\mathfrak A}}
\nc{\fB}{{\mathfrak B}}
\nc{\fC}{{\mathfrak C}}
\nc{\fD}{{\mathfrak D}}
\nc{\fE}{{\mathfrak E}}
\nc{\fF}{{\mathfrak F}}
\nc{\fG}{{\mathfrak G}}
\nc{\fH}{{\mathfrak H}}
\nc{\fI}{{\mathfrak I}}
\nc{\fJ}{{\mathfrak J}}
\nc{\fK}{{\mathfrak K}}
\nc{\fL}{{\mathfrak L}}
\nc{\fM}{{\mathfrak M}}
\nc{\fN}{{\mathfrak N}}
\nc{\fO}{{\mathfrak O}}
\nc{\fP}{{\mathfrak P}}
\nc{\fQ}{{\mathfrak Q}}
\nc{\fR}{{\mathfrak R}}
\nc{\fS}{{\mathfrak S}}
\nc{\fT}{{\mathfrak T}}
\nc{\fU}{{\mathfrak U}}
\nc{\fV}{{\mathfrak V}}
\nc{\fW}{{\mathfrak W}}
\nc{\fZ}{{\mathfrak Z}}
\nc{\fX}{{\mathfrak X}}
\nc{\fY}{{\mathfrak Y}}
\nc{\fa}{{\mathfrak a}}
\nc{\fb}{{\mathfrak b}}
\nc{\fc}{{\mathfrak c}}
\nc{\fd}{{\mathfrak d}}
\nc{\fe}{{\mathfrak e}}
\nc{\ff}{{\mathfrak f}}
\nc{\fg}{{\mathfrak g}}
\nc{\fh}{{\mathfrak h}}
\nc{\fiI}{{\mathfrak i}}  
\nc{\ffi}{{\mathfrak i}}  
\nc{\fj}{{\mathfrak j}}
\nc{\fk}{{\mathfrak k}}
\nc{\fl}{{\mathfrak{l}}}
\nc{\fm}{{\mathfrak m}}
\nc{\fn}{{\mathfrak n}}
\nc{\fo}{{\mathfrak o}}
\nc{\fp}{{\mathfrak p}}
\nc{\fq}{{\mathfrak q}}
\nc{\fr}{{\mathfrak r}}
\nc{\fs}{{\mathfrak s}}
\nc{\ft}{{\mathfrak t}}
\nc{\fu}{{\mathfrak u}}
\nc{\fv}{{\mathfrak v}}
\nc{\fw}{{\mathfrak w}}
\nc{\fz}{{\mathfrak z}}
\nc{\fx}{{\mathfrak x}}
\nc{\fy}{{\mathfrak y}}

\nc{\la}{{\lambda }}

\nc{\La}{{\Lambda }}

\nc{\Coh}{{{\mathcal C}oh}}
\nc{\Loc}{{{\mathcal L}oc}}
\nc{\GR}{{G_\bR}}

\newcommand{\ct}{{T^*\!X}}

\newcommand{\To}{\longrightarrow}

 \providecommand{\curlybrackets}[1]{{\{\!\{#1\}\!\}}}
 \nc{\cb}{\curlybrackets}

\nc{\str}{{\hat\cA}}
\nc{\Spec}{{\on{Spec}}}
\nc{\cext}{{\cE\mathit x\mathit t}}
\nc{\ctor}{{\cT\mathit o\mathit r}}
\nc{\crhom}{{\operatorname{R}\cH\mathit o\mathit m}}
\nc{\chom}{{\cH\mathit o\mathit m}}
\nc{\oh}{{\on{H}}}
\nc{\codim}{{\on{codim}}}
\nc{\Supp}{{\on{Supp}}}
\nc{\coker}{{\on{coker}}}
\nc{\rhD}{{\operatorname{Mod}_{hr}(\cD_X)}}
\nc{\rhDL}{{\operatorname{Mod}_{hr}(\cD_X)_\La}}
\nc{\rhE}{{\operatorname{Mod}_{hr}(\cE_X)}}
\nc{\rhEL}{{\operatorname{Mod}_{hr}(\cE_X)_\La}}
\nc{\rhFE}{{\operatorname{Mod}_{hr}(\hat\cE_X)}}
\nc{\rhFEL}{{\operatorname{Mod}_{hr}(\hat\cE_X)_\La}}
\nc{\op}{{\on{P}}}
\nc{\oM}{{\on{M}}}
\nc{\ba}{\begin{array}}
\nc{\ea}{\end{array}}
\nc{\hs}{\hspace*}
\nc{\cl}{\colon}
\newcommand{\seteq}{\mathbin{:=}}

\newcommand{\scbul}{{\,\raise.4ex\hbox{$\scriptscriptstyle\bullet$}\,}}

\nc{\dT}{{\mathring{T}}{}^*}
\nc{\oY}{{\mathring{Y}}}
\nc{\oLa}{{\mathring{\La}}}
\nc{\eq}{\begin{eqnarray}}
\nc{\eneq}{\end{eqnarray}}
\nc{\Per}{\on{\mathcal{P}\mathit{er}}}

\rc{\setminus}{{-}}

\nc{\bigmid}{\;\mathbin{\rule[-1.8ex]{.5pt}{3.8ex}}\;}

\nc{\cmtk}[1]{\color{red}{{\fbox{K}} #1}\color{black}}

\begin{document}

\title{On the codimension-three conjecture}

\author{Masaki Kashiwara}

\address{Research Institute for Mathematical Sciences, 
Kyoto University,                 
Kyoto, 606--8502, Japan\\
and \\ Department of Mathematical Sciences, Seoul National University,
Seoul, Korea}
\thanks{M.\ Kashiwara was partially supported 
by Grant-in-Aid for Scientific Research (B) 23340005,
Japan Society for the Promotion of Science.}
\email{masaki@kurims.kyoto-u.ac.jp}

\author{Kari Vilonen}

\address{Department of Mathematics, Northwestern University, Evanston, IL
60208, USA \\
and \\ Department of Mathematics, Helsinki University, Helsinki, Finland}
\thanks{K.\ Vilonen was supported by DARPA via AFOSR grant FA9550-08-1-0315}
\email {vilonen@math.northwestern.edu, vilonen@math.helsinki.fi}

\date{July 20, 2010}

\maketitle

\section{Introduction}

In this note we sketch a proof of a fundamental conjecture, 
{\em the codimension-three conjecture}, 
for microdifferential holonomic systems with regular singularities 
This conjecture emerged at the end of the 1970's and is well-known among experts. However, as far as we know, it was never formally written down as a conjecture, perhaps because of lack of concrete evidence for it. As one can also view our result from the point of view of perverse sheaves we have written much of the introduction from  that point of view.

Let $X$ be a complex manifold. 
The notion of micro-supports introduced by \cite{KS1, KS2} allows us
to study perverse sheaves micro-locally 
(i.e., locally on the cotangent bundle).
Let us fix a conic Lagrangian subvariety $\La\subset \ct$. 
It is often important and interesting 
to understand the category $\op_\La(X)$ of perverse sheaves on $X$ 
with micro-support  in $\La$. 
Equivalently, via the Riemann-Hilbert correspondence, 
we can phrase this problem in terms of $\cD_X$-modules. 
From that point of view we can view $\op_\La(X)$ 
as the category of regular holonomic  $\cD_X$-modules 
whose characteristic variety is contained in $\La$. The basic structure of this category has been studied by several authors, for example, 
\cite{Be}, \cite{MV} and \cite{KS2}.
In \cite{GMV1} it is shown 
how in principle one can describe this category: 
it is equivalent to the category of finitely generated modules
over an associative finitely presented algebra. 
However, it is perhaps more interesting to describe $\op_\La(X)$ 
in terms of the geometry of $\ct$. This is the problem we consider here. 

The category $\op_\La(X)$ gives rise to a stack $\Per_\La$ on $\ct$. 
From the point of view of perverse sheaves, it is 
the stack of microlocal perverse sheaves  (see \cite{W}). 
From the point of view of $\cD_X$-modules the construction is more transparent:
one simply passes from regular holonomic $\cD_X$-modules 
to  regular holonomic $\cE_X$-modules; 
here $\cE_X$ stands for the ring of microdifferential operators on $\ct$. 
One expects the microlocal description of $\op_\La(X)$, 
i.e., the description of $\Per_\La(\La)$ to be conceptually simpler. 

Let us write 
\begin{equation}
\La  \ = \ \La^0\sqcup\La^1 \sqcup \La^ 2 \sqcup \cdots\,,
\end{equation}
where $\La^i$ is the locus of codimension $i$ singularities of $\La$. 
We leave the appropriate notion of ``singularity"  vague for now. 
We set $\La^{\geq i}=\cup_{k\ge i}\La^k$.
It is not difficult to show,
either from the topological or
from the analytic point of view, that  the following two statements hold:
\begin{equation}\label{codim1}
\text{The functor $\Per_\La(\La)\To\Per_\La(\La\setminus\La^{\geq 1})$
is faithful,}
\end{equation}
and
\begin{equation}\label{codim2}
\text{The functor $\Per_\La(\La)\to\Per_\La(\La\setminus\La^{\geq 2})$
is fully faithful.}
\end{equation}

In particular the latter implies
\begin{equation}
\parbox{65ex}
{If we have the Lagrangian $\La =\La_1\cup\La_2$ 
with each $\La_i$ Lagrangian and
$\codim_\La(\La_1\cap\La_2 )\geq 2$, then  $\Per_\La(\La) =  
\Per_{\La_1}(\La_1) \times  \Per_{\La_2}(\La_2)$.}
\end{equation}
In concrete terms, \eqref{codim2} means that 
beyond the codimension one singularities of $\La$ 
only conditions on objects are imposed. 
All the essential data are already given along $\La^0$ and $\La^1$. 
Along the locus $\La^0$ we specify a local system (cf.\ \cite{K1}),
and along  $\La^1$  we specify some ``glue" 
between the local systems on various components of  $\La^0$. 
Such a description of the stack $\Per_\La(\La - \La^{\geq 2})$, 
in terms of Picard-Lefschetz/Morse theory is discussed in \cite{GMV2}.

In this paper we answer the question as to what happens 
beyond codimension two, i.e., we announce the following fundamental fact.

\begin{thm}
\label{main theorem}
For an open subset $U$ of $\La$ and a closed analytic subset $Z$ of $U$
of codimension at least $3$,
the functor $\Per_\La(U)\to\Per_\La(U\setminus Z)$ is
an equivalence of categories.
\end{thm}

We sketch an analytic proof of this result utilizing the ring of  microdifferential operators.
The detailed proof  will appear in a forthcoming paper.

The second author wishes to thank Kari Astala, Bo Berndtsson, Laszlo Lempert, Eero Saksman, Bernard Shiffman, Andrei Suslin, and Hans-Olav Tylli for helpful conversations.

\section{The set up}

Let us recall the definition
of  the sheaves of rings of microdifferential operators,
$\cE_X$ and $\hE_X$ (see \cite{SKK, Sch, K2}).
Let $\pi_X\cl T^*X\to X$ be the cotangent bundle to $X$. 
The $\bC^*$-action on $\ct$ gives rise to the Euler vector field $\chi$. 
We say that a function $f(x,\xi)$ on $\ct$ is homogeneous of degree $j$ 
if $\chi f =j f$. Then
we define the sheaf $\hE_X(m)$ for $m\in\Z$ by setting, for an open subset $U$ of $T^*X$, 
\begin{equation*}
\hE_X(m)(U) = \bigr\{\sum\limits_{j=-\infty}^m p_j(x,\xi)\mid
\text{$p_j(x,\xi)\in\cO_{T^*X}(U)$ is homogeneous of degree $j$}\bigr\}
\end{equation*}
and then setting $\hE_X=\bigcup_{m\in\Z}\hE_X(m)$.
The expression $\sum_{j=-\infty}^m p_j(x,\xi)$ 
is to be viewed as a formal expression. 
The formal expressions are multiplied using the Leibniz rule. In this manner
$\hE_X$ becomes a sheaf of rings.

We define $\cE_X$ to be the subsheaf of rings of $\hE_X$ 
consisting of symbols \linebreak $\sum_{j=-\infty}^m p_j(x,\xi)$ 
which satisfy the following growth condition:
\begin{equation}
\begin{gathered}
\text{for every compact $K\subset U$ there exists a $C>0$ such that}
\\
\sum_{j=-\infty}^0 \|p_j(x,\xi)\|_K \frac {C^{-j}}{(-j)!} <\infty\,;
\end{gathered}\label{definition of E}
\end{equation}
here  $\|p_j(x,\xi)\|_K$ stands for the sup norm on $K$. 
Standard estimates can be used to show that $\cE_X$ is 
indeed closed under multiplication 
and hence constitutes a subring of $\hE_X$.

In the study of both $\cE$-modules and $\hE$-modules,
we can make use of canonical transformations. 
In particular, any holonomic module $\cM$ 
can locally be put in general position. We say that $\cM$ is in general position if
its support $\La$ is such that the fibers of the projection $\La \to X$ are at most one dimensional. 
We have the following basic fact:
\begin{thm}
Let us assume that a holonomic $\cE_X$-module $\cM$ 
is in generic position at the point $p\in\ct$. 
Then the local $\cE_{X,p}$-module $\cM_p$, 
the module $\cM$ localized at the point $p$, 
is a holonomic $\cD_{X,\pi_X(p)}$-module,
and the canonical morhism $\cE_{X,p}\otimes_{\cD_{X,\pi_X(p)}}\cM_p\to \cM_p$
is an isomorphism. 
\end{thm}

A proof of this result is given in \cite{Bj}, theorem 8.6.3, for example. 
The proof uses the same reduction that we utilize later 
in this paper combined with Fredholm theory. 
The estimates in the definition \eqref{definition of E} 
are precisely the  ones so that this theorem holds.

Let us recall the notion of regular singularities. 
For a coherent $\cE_X$-module $\cM$,
a coherent $\cE_X(0)$-submodule $\cN$ is called an $\cE_X(0)$-lattice
if $\cE_X\otimes_{\cE_X(0)}\cN\to\cM$ is an isomorphism.
A holonomic $\cE$-module with support $\La$
is said to have {\em regular singularities}
if locally near any point on the support of $\cM$ 
the module $\cM$ has an $\cE(0)$-lattice $\cN$ 
which is invariant under $\cE_\La(1)$, 
the subsheaf of order 1 operators whose principal symbol vanishes on $\La$. 
Kashiwara and Kawai show, using their notion of order:

\begin{thm}
A regular holonomic  $\cE$-module possesses 
a globally defined $\cE(0)$-lattice invariant under  $\cE_\La(1)$. 
The analogous result holds for $\hE$-modules. 
\end{thm}

For a proof see \cite{KK},  Theorem 5.1.6. 
In the rest of the paper we make use of the (global)
existence of an $\cE(0)$-lattice.
Its invariance under $\cE_\La(1)$ will play no role.

Recall that the stack $\Per_\Lambda$ is equivalent to the stack of regular holonomic $\cE_X$-modules with support in $\Lambda$. We prove the extension theorem \ref{main theorem} in the context of regular holonomic $\cE_X$-modules.

\section{Definition of sheaves of rings}

In the next section we reduce the codimension-three conjecture 
to an analogous conjecture about another related ring.
We will define this ring below.  

Let us consider the formal power series ring $\hA = \bC[[t]]$. 
It is, of course, a DVR (discrete valuation ring). 
We will define a subring $A$ of $\hA$ in the following manner. 
For any $C>0$ we define a norm $\|\ \|_C$ on $\hA$ by the formula
\begin{equation}
\|\sum_{j=0}^\infty a_jt^j\|_{_C}\ = \ \sum_{j=0}^\infty |a_j|\frac{C^j}{j!} 
\,.
\end{equation}
We write $A_C$ for the subring consisting of elements of $\hA$ with finite norm. The ring $A_C$ is a Banach DVR. Finally, we set 
\begin{equation}
A = \varinjlim _{C \to 0}A_C\,.
\end{equation}
The ring $A$ is DNF (Dual Nuclear Frechet) DVR. 
We write $K$ for the fraction field of $A$ 
and $\hK$ for the field of Laurent series. 

Let $X$ be a complex manifold. 
We write $\cA_X$ for the sheaf of holomorphic functions on $X$ 
with values in $A$ and similarly, 
for $\widehat \cA_X$, $\cK_X$, and $\widehat\cK_X$. 
We can also view $\cA_X$ as a topological tensor product  
$\cA_X = A \widetilde\otimes_\bC \cO_X$, 
and similarly for  $\widehat \cA_X$, $\cK_X$, and $\widehat\cK_X$.

 \section{The reduction}

 In this section we make a basic reduction of the conjecture.  
In the sequel, we consider the projective cotangent bundle
$P^*X\seteq\dT X/\C^*$ where $\dT X\seteq T^*X\setminus X$.
Since $\cE_X$ and $\hE_X$ are constant along the fibers of
$\dT X\to P^*X$, we regard $\cE_X$ and $\hE_X$ 
as sheaves of rings on $P^*X$.

Let $\La$ be a locally closed Lagrangian subvariety of $P^*X$,
let $\oLa\subset \La$ 
be an open subset such that $\La\setminus\oLa$ is 
a closed analytic subset of codimension at least three in $\La$. 
We assume that we are given a regular holonomic $\cE_X$-module $\cM$ 
with support $\La$ which is defined on  $\oLa$,
and an $\cE_X(0)\vert_{\oLa}$-lattice $\cN$ of $\cM$.
As  the extension problem is local, 
we can work in a neighborhood of a particular point 
$p\in\La\setminus\oLa$ which we might as well assume, 
working by induction,  to be a smooth point of  $\La\setminus\oLa$. 
We now consider the extension problem in the vicinity of this point 
and we put $\La$, via a canonical transformation, 
in generic position at $p$. 
In the neighborhood of $\pi_X(p)$ we take
a coordinate system $(x_1,\ldots, x_n)$
such that $\pi_X(p)$ corresponds to the origin 
and the point $p$ corresponds to $dx_n$ at the origin.  
Then we have a finite map, defined locally in the neighborhood of $p$\;:

\begin{equation}
\rho\cl \La \to Y\subset\bC^{n-1}
\quad
\text{given by $x_1,\ldots, x_{n-1}$.}
\end{equation}
where $Y$ is an open subset of $\bC^{n-1}$.
We may assume that there exists an open subset $\oY$ of $Y$
such that $\rho^{-1}(\oY)=\oLa$ and $Y\setminus\oY$ is an
analytic subset of codimension $\ge3$.
Let us now consider the $\rho_*\cE_X$-module $\rho_*\cM$ 
and the  $\rho_*\cE_X(0)$-module $\rho_*\cN$.  
Also, we write $t$ for the variable $\partial_{x_n}^{-1}$. 
Using this notation $\cA_Y$ is a subsheaf of rings of $\rho_*\cE_X(0)$ 
and  $\cK_Y$ is a subsheaf of rings of $\rho_*\cE_X$ and similarly in the formal case. By a standard argument using division theorems (see \cite{SKK}),
we conclude that

\begin{prop} 
The sheaf  $\rho_*\cM\vert_{\oY}$ is coherent over  $\cK_Y\vert_{\oY}$ 
and  $\rho_*\cN\vert_{\oY}$ is coherent over $\cA_Y\vert_{\oY}$.
\end{prop}

Let $i\cl \oLa\hookrightarrow\La$ and
$j\cl \oY\hookrightarrow Y$ be the open inclusions.
Then a standard argument shows 
\begin{equation}
\parbox{65ex}{The $\cE_X$-module $i_*\cM$ is a coherent $\cE_X$-module
if and only if 
 $j_*\rho_*\cM$ is a coherent $\cK_Y$-module.}
\end{equation}

 We now recall that holonomic modules are Cohen-Macaulay, i.e., we have 
\begin{equation}
\cext^k_{\cE_X}(\cM,\cE_X)\vert_{\oLa} = 0 \qquad \text{unless} \ \ k=\dim X\,.
\end{equation}

This immediately implies that $\rho_*\cM$ is flat as a $\cK_Y$-module:
\begin{equation}\label{eq:van}
\cext^k_{\cK_Y}(\cM,\cK_Y)\vert_{\oY} = 0 \qquad \text{unless} \ \ k=0\,.
\end{equation}

In particular, we conclude
\begin{equation}\label{eq:free}
\text{$\cM\vert_{\oY}$ is locally free as a $\cK_Y\vert_{\oY}$-module}\,.
\end{equation}

Note that statement \eqref{eq:van}
does not imply \eqref{eq:free} immediately. 
However, we can conclude \eqref{eq:free} by using a  result 
of Popescu, Bhatwadekar, and Rao \cite{P}; for a nice discussion, see also \cite{S}. They show:
\begin{equation}
\parbox{70ex}{Let $R$ be a regular local ring containing a field 
with maximal ideal $\fm$ and $t\in\fm\setminus\fm^2$. 
Then every finitely generated projective module 
over the localized ring $R_t$ is free.}
\end{equation}
This result is related to Serre's conjecture and was conjectured by Quillen in \cite{Q}.

\section{The formal case}

In this section we discuss the proof of the theorem in the formal case. 
In the previous section we reduced the proof to the following situation. Let $Y$ be a complex manifold and $\oY$ an open subset of $Y$ 
such that $Y\setminus\oY$ is an analytic subset
of codimension~$\ge3$.
We write $j\cl \oY\to Y$ for the open inclusion. 
We assume that we are given a locally free $\widehat\cK_Y\vert_{\oY}$-module 
$\hM$ of finite rank
together with an $\widehat\cA_Y\vert_\oY$-lattice 
$\hN$ in  $\hM$. Let us write 

\begin{equation}
\hN^* \ =  \ \chom_{\cA_\oY}(\hN,\cA_\oY)\,.
\end{equation}
By replacing $\hN$ with $\hN^{**}$, 
we may assume that $\hN$ is reflexive
(i.e., that $\hN\to\hN^{**}$ is an isomorphism). 

We claim:

\begin{thm}
If $\hN$ is a reflexive coherent $\str_{\oY}$-module,
then $j_*\hN$ is a coherent $\str_Y$-module. 
\end{thm}

From this theorem we immediately conclude 
that $j_*\hM$ is a coherent $\widehat\cK_Y$-module 
and hence obtain our theorem in the formal case. 

As to the proof of the theorem above, we reduce it, using \cite[Proposition 1.2.18]{KS3},
to the following
analogous classical result in several complex variables. 
We state this result, due to Trautman, Frisch-Guenot, 
and Siu, in a special case suitable for us. 
\begin{thm}
If $\cF$ is a reflexive coherent $\cO_{\oY}$-module, 
then $j_*\cF$ is a coherent $\cO_Y$-module. 
\end{thm}

For a discussion of results of this type, 
see the S\'eminaire Bourbaki talk by Douady (\cite{D}). 

\section{Convergent case}

In this section we will argue that formal case 
implies the convergent case. 
Just as in the previous section we are reduced to the following situation. 
Let $Y$ be a complex manifold and $\oY$ an open subset of $Y$ 
such that $Y\setminus \oY$ is an analytic subset 
of codimension $\ge3$. 
We write $j\cl\oY\to Y$ for the open inclusion. 
We assume the we are given a locally free $\cK_{\oY}$-module $\cM$ 
of finite rank on $\oY$. 
We assume that it is locally free as a 
$\cK_Y$-module.
As the corresponding $\hE_X$-module extends across $\La\setminus\oLa$, 
we obtain the additional piece of information 
that the sheaf $j_*\hM$ is locally free.
As we work locally along $Y\setminus \oY$, we can then assume 
that $j_*\hM$ is free.
Thus, we have, locally along $Y$ the following situation:
\begin{equation}
\begin{gathered}
\xymatrix{
 {\hM \cong \widehat\cK_{\oY}^{\oplus r}}\ar@{}[dr]|{\square}
& \widehat\cA_{\oY}^{\oplus r} 
\ar@{_{(}->}[l] \\
 \cM\ar@{_{(}->}[u] &
{\strut\;(\cM\cap  \widehat\cA_{\oY}^{\oplus r}) =\cN}\ar@{_{(}->}[l]
\ar@{_{(}->}[u]
}
\end{gathered}
\end{equation}
As $\cM$ is locally free on $\oY$ 
we conclude that the  $\cK_{\oY}$-lattice $\cN$ 
is also locally free on  $\oY$. We will now conclude 
using the following result.
\begin{thm}
\label{formal comparison}
Let Y be a complex manifold and $\oY$ 
an open subset of $Y$ such that
$Y\setminus\oY$ is a closed submanifold of codimension $\ge 2$.
Then any locally free $\cA_{\oY}$-module $\cN$ 
whose formal completion $\hN$ extends as a locally free module on $Y$ 
also extends to a locally free module on $Y$. 
\end{thm}

In the proof of this theorem we make use of the following fact: 
\eq&&
\parbox{68ex}{Let $B$ be a Banach algebra and let
$\cO^B_X$ (resp.\ $\cC_X^B$) 
be the sheaf of rings of $B$-valued holomorphic (resp.\ continuous) 
functions on
a Stein manifold $X$. Then two 
locally free $\cO^B_X$-modules of finite rank are 
isomorphic as soon as they are isomorphic after tensoring by
$\cC_X^B$.}
\eneq
 This an extension of the Oka principle to infinite dimensional bundles. These issues are discussed in the paper \cite{Bu}. 

The following result on extensions of submodules also 
follows from theorem \ref{formal comparison}:
\begin{thm} 
Let $\cM$ be a holonomic $\cE_X$-module defined on
an open set $\Omega$ of $T^*X$ and let $Z$ be a closed analytic subset
of $\Omega$ of codimension at least $\dim X+2$.
Then any holonomic submodule $\cN$ of $\cM\vert_{\Omega\setminus Z}$
extends to a holonomic $\cE_X\vert_\Omega$-submodule of $\cM$.
\end{thm} 
Note that if an extension exists it is unique by 
\eqref{codim2}.

 \section{Open problems}

 In this section we discuss open problems which are closely related to our main result. 
 
 We have the following result:
 \begin{prop}
The category of regular holonomic $\cE_X$-modules is 
a full subcategory of the  category of regular holonomic $\hE_X$-modules
\end{prop}

Thus, it is natural to conjecture:

\begin{conj}
The category of regular holonomic $\cE_X$-modules is equivalent 
to the  category of regular holonomic $\hE_X$-modules
\end{conj}

If we fix supports and consider the subcategories 
where the objects have a fixed conic Lagrangian support $\La$,
one can check that the two categories 
coincide in codimension zero by direct verification.
Then, by theorem \ref{formal comparison}, 
the problem reduces to the case of pure codimension one locus. 

\medskip

Microlocal perverse sheaves can be defined with arbitrary coefficients, 
for example, with integral coefficients. 
This gives us a stack $\Per_\La(\bZ)$ on a Lagrangian $\La$. We conclude with:

\begin{conj}
The codimension-three conjecture holds for the stack $\Per_\La(\bZ)$.
\end{conj}


\begin{thebibliography}{GMV2}
 
\bibitem[Be]{Be}

Alexander Beilinson,
{\em How to glue perverse sheaves},
K-theory, arithmetic, and geometry (Moscow 1984), Lecture notes in Math., vol 1289, Springer, Berlin, 1987, 42--51.

\bibitem[Bj]{Bj}
Jan-Eric Bj\"ork,
{\em Analytic $\cD$-modules and applications},
Mathematics and its Applications, 247.
Kluwer Academic Publishers Group, 1993. 581 pp.


 \bibitem[Bu]{Bu}  Lutz Bungart,
{\em On analytic fiber bundles I.
\ Holomorphic fiber bundles with infinite dimensional fibers},
Topology {\bf7} 1967 55--68.
 
 \bibitem[D]{D}Adrien Douady, 
{\em Prolongement de faisceaux analytiques coherents 
(travaux de Trautmann, Frisch-Guenot et Siu)}, S\'eminaire Bourbaki, 1969/70.

\bibitem[GMV1]{GMV1} Sergei Gelfand, Robert MacPherson and Kari Vilonen,
{\em  Perverse sheaves and quivers},
Duke Math. J. {\bf 83} (1996), no. 3, 621--643.


\bibitem[GMV2]{GMV2} 
\bysame,
{\em Micro-local perverse sheaves},
arXiv:math.AG/0509440.

\bibitem[K1]{K1} Masaki Kashiwara, 
{\em Introduction to Microlocal Analysis},
l'enseignement math\'ematique {\bf 32}
(1986) 227--259.

\bibitem[K2]{K2}
\bysame,
{\em $D$-modules and microlocal calculus},
Translated from the 2000 Japanese original by Mutsumi Saito,
Translations of Mathematical Monographs, {\bf217},
Iwanami Series in Modern Mathematics,
American Mathematical Society, Providence, RI (2003).

\bibitem[KK]{KK}  Masaki Kashiwara and Takahiro Kawai,
{\em On holonomic systems of microdifferential equations. III,
Systems with regular singularities},
Publ. Res. Inst. Math. Sci. {\bf 17} (1981), no. 3, 813--979.

\bibitem[KS1]{KS1}
Masaki Kashiwara and Pierre Schapira,
{\em Micro-support des faisceaux: application aux modules diff\'erentiels},
C. R. Acad. Sci. Paris S\'er. I Math.  {\bf295}  (1982), no. 8, 487--490.

\bibitem[KS2]{KS2} \bysame,
{\em Sheaves on Manifolds}, Grundlehren der Mathematischen 
Wissenschaften {\bf 292}, Springer-Verlag
(1990).

\bibitem[KS3]{KS3}
\bysame,
{\em Deformation quantization modules},  
arXiv:1003.3304.
											     
\bibitem[MV]{MV} Robert MacPherson and Kari Vilonen,
{\em Elementary construction of perverse sheaves}, 
Invent. Math. {\bf 84} (1986), no. 2, 403--435. 

\bibitem[P]{P} Dorin Popescu,
{\em On a question of Quillen},
Bull. Math. Soc. Sci. Math. Roumanie (N.S.) 45(93) (2002), no. 3-4, 
(2003) 209--212 .

\bibitem[Q]{Q} Daniel Quillen, 
{\em  Projective modules over polynomial rings},
 Invent. Math. 36 (1976), 167--171. 

\bibitem[SKK]{SKK}
Mikio Sato, Takahiro Kawai and Masaki Kashiwara,
{\em Microfunctions and pseudo-differential equations},
Hyperfunctions and pseudo-differential equations 
(Proc. Conf., Katata, 1971; dedicated to the memory of Andr\'e Martineau),  
Lecture Notes in Math., {\bf287}, Springer, Berlin, (1973)
265--529. 


\bibitem[Sch]{Sch}
Pierre Schapira, 
{\em Microdifferential Systems in the Complex Domain},
Grundlehren der Mathematischen Wissenschaften 
{\bf269} Springer-Verlag, Berlin (1985).


\bibitem[S]{S} Richard Swan, 
{\em N\'eron-Popescu desingularization},
Algebra and geometry (Taipei, 1995), 
Lect. Algebra Geom., 2, Int. Press, Cambridge, MA, (1998) 135--192.

\bibitem[W]{W} Ingo Waschkies, 
{\em The stack of microlocal perverse sheaves},
Bull. Soc. Math. France {\bf 132}  (2004),  no. 3, 397--462. 

\end{thebibliography}
\end{document}